\documentclass[12pt]{amsart}

\def\@typesizes{%
       \or{5}{6.5}\or{6}{7.5}\or{7}{8.5}\or{8}{11}\or{9}{12}%
       \or{10}{13}
       \or{\@xipt}{14}\or{\@xiipt}{15}\or{\@xivpt}{18}%
       \or{\@xviipt}{20}\or{\@xxpt}{24}}

\usepackage{amsthm}
\usepackage{amsmath}
\usepackage{amssymb}
\usepackage{graphicx}
\usepackage{enumerate}
\usepackage{color}
\allowdisplaybreaks
\numberwithin{equation}{section}
\numberwithin{figure}{section}

\theoremstyle{plain}
\newtheorem{theorem}{ Theorem}[section]
\newtheorem{proposition}[theorem]{ Proposition}
\newtheorem{lemma}[theorem]{ Lemma}
\newtheorem{corollary}[theorem]{ Corollary}
\newtheorem{example}[theorem]{ Example}
\newtheorem{remark}[theorem]{ Remark}
\newtheorem{definition}[theorem]{ Definition}
\newtheorem{conjecture}{ Conjecture}

\def\BET{\begin{theorem}}
\def\ENT{\end{theorem}}
\def\BEP{\begin{proposition}}
\def\ENP{\end{proposition}}
\def\BEL{\begin{lemma}}
\def\ENL{\end{lemma}}
\def\BEC{\begin{corollary}}
\def\ENC{\end{corollary}}
\def\BEE{\begin{example} \rm}
\def\ENE{\end{example}}
\def\BER{\begin{remark} \rm}
\def\ENR{\end{remark}}
\def\BED{\begin{definition} \rm}
\def\END{\end{definition}}
\def\BECJ{\begin{conjecture}}
\def\ENCJ{\end{conjecture}}

\def\bea{\begin{eqnarray}}
\def\eea{\end{eqnarray}}

\def\beas{\begin{eqnarray*}}
\def\eeas{\end{eqnarray*}}

\def\beq{\begin{equation}}
\def\eeq{\end{equation}}

\def\beal{\begin{align*}}

\def\eeal{ \end{align*} }

\def\row{\nonumber \\ & & }
\def\roweq{\nonumber \\ &=& }
\def\rowleq{\nonumber \\  & \leq & }

\def\rowpl{\nonumber \\ & + & }
\def\rowmi{\nonumber \\ & - & }

\def\bbC{{\mathbb C}}

\def\bbN{{\mathbb N}}

\def\bbR{{\mathbb R}}
\def\bbS{{\mathbb S}}

\def\cH{{\mathcal H}}

\def\cO{{\mathcal O}}

\def\cS{{\mathcal S}}
\def\cT{{\mathcal T}}

\def\cV{{\mathcal V}}

\def\cY{{\mathcal Y}}
\def\cZ{{\mathcal Z}}

\def\ov{\overline}
\def\ef{\eqref}

\begin{document}

\title[Blinking Steklov eigenvalues]{"Blinking eigenvalues" of the Steklov 
problem  generate the continuous spectrum in a cuspidal domain}



\author{Sergei A. Nazarov}
\address{
Saint-Petersburg State University,
Universitetskaya nab., 7--9,  St. Petersburg, 199034, Russia, and  \hfill\break
Institute of Problems of Mechanical Engineering RAS,
V.O., Bolshoj pr., 61, St. Petersburg, 199178, Russia}
\email{s.nazarov@spbu.ru, srgnazarov@yahoo.co.uk}

\author{Jari Taskinen}
\address{Department of Mathematics and Statistics, P.O.Box 68, 
University of Helsinki, 00014 Helsinki, Finland}
\email{jari.taskinen@helsinki.fi}

\thanks{ The first  named author was  supported by the grant 
17-11-01003 of the Russian Science Foundation.} 

\begin{abstract}
We study the Steklov spectral problem for the Laplace operator 
in a bounded domain $\Omega \subset \bbR^d$, $d \geq 2$, with a cusp
such that the  continuous spectrum of the problem is non-empty,
and also in the family of bounded domains $\Omega^\varepsilon \subset \Omega$, 
$\varepsilon > 0$, obtained from $\Omega$ by  blunting the cusp at the 
distance of $\varepsilon$ from the cusp tip. 
While the spectrum in the blunted domain $\Omega^\varepsilon$  consists 
for a fixed $\varepsilon$  of  an unbounded positive sequence 
$\{ \lambda_j^\varepsilon \}_{j=1}^\infty$ of eigenvalues, 
we single out different types of behavior of some 
eigenvalues  as  $\varepsilon \to  +0$: in particular, stable, blinking, 
and gliding
families of eigenvalues are found.  We also describe  a mechanism which
transforms the family of the eigenvalue sequences into the continuous 
spectrum of the problem in $\Omega$, when $\varepsilon \to +0$.

\end{abstract}

\maketitle


\section{Introduction.}
\label{sec1}
\subsection{Formulation of the problems.}
\label{sec1.1}

\begin{figure}
\label{fig1}
\begin{center}
{\includegraphics[height=6cm,width=10cm]{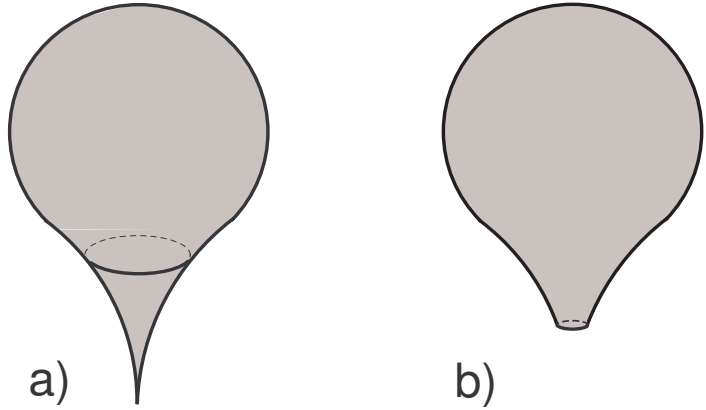}} 
\caption{a) Cuspidal domain, b) domain with blunted cusp. }
\end{center}
\end{figure}

Let $\Omega$ be a domain  in $\bbR^{n}$, $n \geq 2$, with  compact closure
$\overline{\Omega}$ and  boundary $ \partial \Omega$ which is smooth 
everywhere 
except at the origin $\cO$ of the Cartesian coordinate system $x=(y,z) \in \bbR^{n-1} \times \bbR$ (Fig.\,\ref{fig1},a). In a neighborhood of the point $\cO$ the domain $\Omega$ coincides with 
the cusp
\bea
\Pi^d = \{ x = (y,z)  \, : \, z \in (0,d) 
, \eta = z^{-m }y \in \omega\} \ ,  \ \ d > 0 , \label{0}
\eea
where  $m > 1$ is the  sharpness exponent of the cusp 
and the cross-section $\omega$ is a domain 
in $\bbR^{n-1}$ bounded by a smooth $(n-2)$-dimensional closed surface
$\partial \omega$. 

First of all, we consider the Steklov problem for the Laplace operator
\bea
- \Delta u (x)&=& \lambda u(x), \, x \in \Omega , \ \ \ \partial_\nu u (x)  = \lambda  u(x),  \, 
x \in \partial\Omega ,
\label{1}
\eea
where $\partial_\nu$ is the outward normal derivative and $\lambda$ is the 
spectral parameter. 

We introduce the Hilbert space $\cH$ endowed with the norm
\bea
\Vert u ; \cH \Vert = \big( \Vert \nabla u ; L^2(\Omega) \Vert^2 + 
\Vert  u ; L^2(\partial \Omega) \Vert^2 \big)^{1/2} \label{2}
\eea
and contained in the Sobolev space $H^1(\Omega)$. Then, the integral identity 
corresponding to the problem \eqref{1} reads as
\bea
( \nabla u, \nabla v)_\Omega = \lambda (u,v)_\Omega \hskip1cm \forall \, 
v \in \cH , \label{3}
\eea
see  \cite{Lad}. Here, $\nabla $ is the gradient, $( \ \,  , \ )_\Upsilon$ is the natural 
scalar product in the Lebesgue space $L^2(\Upsilon)$, while the scalar product in
$\cH$ generated by the norm \eqref{2} will be denoted by $\langle \ \, , \ \rangle$ in the 
following. Moreover, we define the operator $\cS : \cH \to \cH$ and the new
spectral  parameter $\mu$ by
\bea
\langle \cS u,v\rangle &=& (u,v)_{\partial \Omega } \hskip1cm \forall \, u,v \in \cH , 
\label{4} \\
\mu &=& (1+ \lambda)^{-1} , \label{5}
\eea
and by these relations the problem \eqref{3} is converted to the abstract equation
\bea
\cS u = \mu u \hskip1cm \mbox{in} \ \cH. \label{6}
\eea
Clearly, the operator $\cS$ is positive definite, symmetric and continuous, and, 
therefore, self-adjoint. 

If we assume for a while that $m \leq1$ in \eqref{0}, the boundary $\partial \Omega$
becomes Lipschitz and the essential spectrum of $\cS$ consists only of the single
point $\mu = 0$ due to  the compactness of the embedding $\cH = H^1(\Omega) 
\subset  L^2(\Omega)$, cf. \cite[Thm.\,10.1.5.]{BiSo}. The remaining part of the
spectrum is discrete and forms a positive sequence converging to zero so that
according to the relation \eqref{5} the whole spectrum $\sigma$
of the problem \eqref{3} (the Steklov problem \eqref{1}) consists of an
unbounded positive sequence of normal eigenvalues. As verified in \cite{na401}, the spectrum remains discrete, if $m <2$.

However, in the case $m \geq 2$ 
the above-mentioned embedding $ \cH \subset L^2(\Omega)$ loses its compactness, 
see e.g. \cite{na401} and 
\cite{MazPob}, and hence the continuous  components $\sigma_{\rm co}$ of the 
spectra of the operator $\cS$ and the Steklov problem become non-empty. 
The component $\sigma_{\rm co} = [\lambda_\dagger , + \infty )$ will be described 
explicitly in Section \ref{sec4} for the most interesting case 
\bea
m=2 ,   \label{7}
\eea
where the positive cut-off value $\lambda_\dagger$ will be obtained from 
\eqref{27}.
Note that in the case $m >2$ it was shown in \cite{na401} that $\lambda_\dagger =0$
and $\sigma_{\rm co} = [0, + \infty)$; this case will not be considered in the 
present paper.

Blunting the cuspidal tip makes the boundary Lipschitz again, see Fig.\,\ref{fig1},b. We consider the simplest truncation
\bea
\Omega^\varepsilon = \Omega \smallsetminus \overline{\Pi^\varepsilon}
\hskip1cm \mbox{with a small} \ \varepsilon \in (0,d) \label{7A}
\eea
and the mixed boundary-value problem
\bea
- \Delta u^\varepsilon (x) &=& 0  , \hskip1.5cm x \in \Omega^\varepsilon, \label{8}  \\
\partial_\nu u^\varepsilon (x)&=& \lambda^\varepsilon  u^\varepsilon (x) ,
\hskip0.4cm x \in ( \partial \Omega)^\varepsilon = \partial \Omega^\varepsilon
\smallsetminus \overline{\omega^\varepsilon}, \label{9}  \\
u^\varepsilon (x) &=& 0 , \hskip1.5cm x \in \omega^\varepsilon, \label{10}  
\eea
with the artificial Dirichlet condition in the end $\omega^\varepsilon = 
\{x \in \Pi^d \, : \, z = \varepsilon\}$. Other truncation surfaces and
types of the artificial boundary condition will be discussed in Section \ref{sec4}.

The operator formulation of the problem \eqref{8}--\eqref{10} will be given in 
Section \ref{sec4}, but is is clear that its spectrum $\sigma^\varepsilon$ is 
discrete and consists of the positive unbounded sequence of eigenvalues 
\bea
0 < \lambda_1^\varepsilon < \lambda_2^\varepsilon 
\leq \lambda_3^\varepsilon \leq \ldots \leq \lambda_m^\varepsilon \leq \ldots
\to + \infty . \label{11}
\eea

The main goal of our paper is to describe the abnormal behavior of some 
entries in 
\eqref{11}, when $\varepsilon \to +0$ and the domain sharpens into a cusp. 
Furthermore, we will find  a mechanism transforming the family of the sequences 
\eqref{11} into the continuous spectrum $\sigma$ of the original Steklov problem 
\eqref{1}.

We will not investigate the asymptotics of all eigenvalues \eqref{11} but 
only some  of them. First,  in Section \ref{sec4.3} 
we find families of eigenvalues which  have the property that
for some $\lambda^{\sf tr}$ and {\it any} small enough  $\varepsilon > 0$, 
the $c \varepsilon$-neighborhood of $\lambda^{\sf tr}$ contains an eigenvalue
belonging to the sequence \eqref{11}, for some positive constant
$c$ independent of $\varepsilon$. For brevity, we call such 
families "stable eigenvalues". (In Section \ref{sec4.1},$3^\circ$, we
make a remark showing that every $\lambda > \lambda_\dagger$ indeed 
is an eigenvalue of the problem \eqref{8}--\eqref{10} for some $\varepsilon$.)

Moreover, Theorem \ref{ASY2} shows that any point 
$\lambda > \lambda_\dagger$ becomes a  "blinking eigenvalue" (Section 
\ref{sec4.1},$2^\circ$) when $\varepsilon \to +0$, i.e. there exists a 
positive  sequence $\{ \varepsilon_k \}_{k=1}^\infty= \{ \varepsilon_k 
(\lambda) \}_{k=1}^\infty$ tending to 0 such that for $\varepsilon = 
\varepsilon_k$, the 
$c_\lambda \varepsilon_k$-neighborhood of $\lambda$ contains an 
entry $\lambda_{m_k}^{\varepsilon_k}$, where $m_k =m_k(\lambda)$. 
(However, for $\varepsilon \not= \varepsilon_k$, there is no guarantee of this
family staying near $\lambda$. The number $\lambda$ becomes a true
eigenvalue of the problem \eqref{8}--\eqref{10} for some $\varepsilon$
close to any entry of the sequence $\{ \varepsilon_k \}_{k=1}^\infty$). 
This fact can 
obviously be used  for the construction of a singular Weyl sequence for the 
operator $\cS$ at the point \eqref{5} (Section \ref{sec5.1}). It is a remarkable 
fact that the structure of the elements of this singular 
sequence is quite  different from the one used in \cite{na401} for the 
continuous spectrum $\sigma_{\rm co}$. 

One more strange phenomenon on the behavior of the eigenvalues of the problem 
\eqref{8}--\eqref{10} will be described in Section \ref{sec4.1},$3^\circ$, 
namely so called "gliding eigenvalues".  We will detect a set  of eigenvalues
$\lambda_{m_k(\varepsilon)}^\varepsilon$, with changing numbers
$m_k(\varepsilon)$, falling down at a high speed $O\big(\varepsilon^{-1} |\ln 
\varepsilon |^{-1} (\lambda_{m(\varepsilon)}^\varepsilon - \lambda_\dagger)
\big) $ as $\varepsilon \to +0$. The speed however declines while approaching
the threshold, which produces a smooth touchdown of 
$\lambda_{m(\varepsilon)}^\varepsilon$ at $\lambda_\dagger$. 
Furthermore, these eigenvalues "sweep" the semi-axis
$(\lambda_\dagger , + \infty)$ many times, when 
$\varepsilon \to +0$ and the Lipschitz domain $\Omega^\varepsilon$ becomes 
cuspidal. (Notice that the "blinking" and "gliding" behaviors do not 
constitute a classification or define separate values of $\varepsilon$ or 
$\lambda$ --- they are just different aspects among the families of eigenvalues.)

The number $\lambda_\dagger$ still belongs to the continuous spectrum, 
since,  according to the general results in \cite{Ko,MaPl2},
see also \cite[Ch.\,10]{KoMaRo1}, eigenvalues of infinite multiplicity do 
not appear in elliptic problems in cuspidal domains so that the essential 
and continuous spectra coincide and thus the latter is also a closed set.

\section{Known facts.}
\label{sec2}
\subsection{Formal asymptotic procedure.}
\label{sec2.1} 
For an eigenfunction of the problem \eqref{1}, we introduce the standard 
asymptotic  ansatz in the analysis of thin domains, which has 
in particular been used in  \cite{na401,na549} 
\bea
u(y,z) = w(z) + W( \eta, z) + \ldots , \label{21}
\eea
where $w$ and $W$ are the power-law functions
\bea
w(z) = z^\tau w_0\ , \ \ W(\eta, z) = z^{\tau +2 } W_0(\eta) , \label{22}
\eea
$\eta = (\eta_1, \ldots ,\eta_{n-1}) = z^{-2} y $ are the stretched coordinates 
in \eqref{0}, and the dots stand for higher-order terms to be  estimated in 
Section \ref{sec4}. We insert \eqref {21}, \eqref{22} to the restriction of the 
problem \eqref{1} on the cusp \eqref{0} and
collect the terms of order $z^{\tau -2}$ in the Laplace equation, and thus 
obtain the $(n-1)$-dimensional Poisson equation with the parameter $z>0$,
\bea
- \Delta_\eta W_0(\eta) = F(\eta) := z^{2 -\tau}\partial_z^2 w(z) \ , \ \ \eta \in 
\omega .
\label{23}
\eea
The unit outward normal vector on the lateral side 
$\Gamma^d = \{ x \, : \, \eta \in \partial \omega , z \in (0,d) \}$
of the cusp equals
\bea
\nu(y,z) = \big( 1 + 4 z^2 |\eta \cdot \nu'(\eta) |^2 \big)^{-1/2}
\big( \nu'(\eta), - 2 z \eta \cdot \nu'(\eta) \big) , \label{23A}
\eea
where $\nu' =( \nu_1' , \ldots, \nu_{n-1}')$ is 
the normal on the boundary $\partial \omega \subset \bbR^{n-1}$ and $|\nu'|=1$. 
Thus, extracting terms of order $z^\tau$ from the Steklov condition yields
the boundary condition 
\bea
\partial_{\nu'} W_0(\eta) = G(\eta) :=
z^{-\tau} \big( \lambda  w(z) + 2 z \eta \cdot \nu'(\eta)  \partial_z
w(z)\big) \ , \ \ \eta \in \partial \omega , \label{24}
\eea
where $\partial_z = \partial / \partial z$, $\partial_{\nu'} = \nu' \cdot 
\nabla_\eta$ and the central dot stands for the scalar product in the Euclidean 
space. According to \eqref{22}, the right-hand sides of \eqref{23} and \eqref{24}
are indeed independent of $z$. The combatibility 
condition in the Neumann problem \eqref{23}, \eqref{24}  is written as
\bea
0 &=& \int\limits_\omega F(\eta) d\eta + \int\limits_{\partial \omega} G(\eta) 
ds_\eta 
\roweq
z^{2 -\tau} \big(|\omega|  \partial_z^2 w(z) + \lambda |\partial  \omega| w(z) + 
2(n-1)  |\omega|^{-1} z \partial_z w(z) \big) \label{24A}
\eea
with the volume $|\omega| = {\rm mes}_{n-1} \omega$ and the area $|\partial \omega|
= {\rm mes}_{n-2} \partial \omega$. We multiply \eqref{24A} by $z^{2n + \tau}$
and, as a result, obtain the ordinary differential 
equation of Euler type
\bea
- \partial_z \big( z^{2(n-1)} |\omega| \partial_z w(z) \big) = 
\Lambda z^{2(n-2)} |\partial \omega| w(z) \ , \ \ z > 0, \label{25}
\eea
with the coefficient
\bea
\Lambda = \frac{|\partial \omega|}{|\omega|} \lambda . \label{25A}
\eea
It has the solutions
\bea
w_\pm (z) = w_0 z^{\tau_\pm} \ \ \mbox{with} \ \tau_\pm 
= - \Big( n - \frac32\Big) \pm
\sqrt{\Big(n-\frac32\Big)^2 - \Lambda } . \label{26}
\eea
The imaginary parts of both exponents $\tau_\pm$ are nonzero provided 
\bea
\Lambda > \big( n - \frac32 \Big)^2 =: \Lambda_\dagger 
\ \ \Leftrightarrow \ \ \lambda > \lambda_\dagger := \Big( n -\frac32\Big)^2
\frac{|\omega|}{|\partial \omega|} , \label{27}
\eea
but both $\tau_\pm$ are real in the case $\lambda < \lambda_\dagger$. Finally, for 
$\lambda=\lambda_\dagger$, the general solution of \eqref{25} is
\beas
w(z) = z^{-n+3/2} ( c_0 + c_1 \ln z). 
\eeas
In Section \ref{sec3} it will be convenient to set
\bea
w_\pm (z) = w_0 z^{-n+3/2} (1 \pm i \ln z) \ \ \ \mbox{at} \ \lambda
= \lambda_\dagger \label{28A}
\eea
so that we can write the general solution for every $\lambda \geq \lambda_\dagger$
as 
\bea
w(z) = b_+ w_+(z) +  b_- w_-(z) \ \ \ \mbox{with} \ b_\pm \in \bbC. \label{28}
\eea
The normalization factor $w_0$ of \eqref{26} and \eqref{28}  will be fixed in 
the formulas \eqref{44}, \eqref{46} of Section \ref{sec3.1}.

Since the compatibility condition \eqref{24A} is satisfied by the function
\eqref{28}, the Neumann problem \eqref{23}, \eqref{24} has a solution $W$
which  is defined up to an additive constant and becomes unique by requiring the 
orthogonality condition 
\bea
\int\limits_\omega W_0(\eta) d\eta = 0 . \label{29A}
\eea

\subsection{The spectrum of the Steklov problem.}
\label{sec2.2}
The following result was proven in \cite{na401} by constructing Weyl singular 
sequences for $\lambda > \lambda_\dagger$ and parametrices for the Steklov problem 
operator  in the case  $\lambda \in [0, \lambda_\dagger)$. 

\BET
\label{T1}
The continuous spectrum $\sigma_{\rm co}$ of the Steklov problem \eqref{1} in the 
cuspidal domain $\Omega$ with the sharpness exponent
\eqref{7} equals  $ [\lambda_\dagger , +
\infty)$, where the cut-off point $\lambda_\dagger$ is given in \eqref{27}. 
\ENT

In other words, the essential spectrum of the operator $\cS$ consists of the point
$\mu=0$ and the continuous spectrum $(0, \mu_\dagger]$,
where $\mu_\dagger =(1 + \lambda_\dagger)^{-1}$, according to \eqref{5}.

Null is an eigenvalue of the problem \eqref{1}, and the interval $(0,
\lambda_\dagger)$  below the continuous spectrum $\sigma_{\rm co}$ may contain 
other points of the discrete spectrum $\sigma_{\rm di}$. Furthermore, it was 
verified in \cite{na401} that in the mirror symmetric case 
\bea
\Omega = \{ (y,z) \, : \,   (-y_1, y_2, \ldots , y_{n-1} , z \} \in \Omega \} \label{29}
\eea
the point spectrum $\sigma_{\rm po}$ is non-empty and in particular it includes the 
unbounded monotone sequence 
\bea
0 < \lambda_1^+ < \lambda_2^+ \leq \lambda_3^+ \leq \ldots \leq \lambda_p^+ \leq 
\ldots \to + \infty \label{30}
\eea
of eigenvalues of the auxiliary problem
\beas
- \Delta u^+ (x) &=& 0, \, x \in \Omega^+ , \ \ \ \partial_\nu u^+ (x) = 
\lambda^+ u^+(x),  \,   x \in (\partial \Omega)^+ ,  \\
 u^+ (x) &=& 0, \,   x \in \Sigma = \partial \Omega^+ \cap \Omega , 
\eeas
where $\Omega^+ = \{ x \in \Omega \, : \, y_ 1 > 0 \}$ is the half-domain and
$\Sigma $ is the artificial truncation surface.

\subsection{Weak formulation of the inhomogeneous Steklov problem.}
\label{sec2.3}
We fix the value of the parameter $\lambda$ in the problem
\bea
- \Delta u(x) = f(x), \, x \in \Omega, \ \ \  
\partial_\nu u(x) - \lambda u(x) = g(x), \,  x \in\partial \Omega
\smallsetminus \cO, \label{31}
\eea
and introduce the space $V_\beta^1(\Omega)$ with the weighted norm
\bea
\Vert u ; V_\beta^1(\Omega) \Vert =
\big( \Vert \nabla u ; L_\beta^2(\Omega) \Vert^2 + 
\Vert u ; L_{\beta-1}^2 (\Omega) \Vert^2 +
\Vert u ; L_\beta^2(\partial \Omega) \Vert^2 \big)^{1/2} \label{32A}
\eea
where $\beta \in \bbR$, $L_\beta(\Upsilon)$ is the weighted Lebesgue space
with the norm
\beas
\Vert v ; L_\beta(\Upsilon)\Vert = \Vert r^\beta v ; L(\Upsilon)\Vert
\eeas
and $r = |x|$ is the distance of $x$ from the cusp tip $\cO$.

\BEL
\label{L1} For every $u \in C_c^\infty(\ov \Omega \smallsetminus \cO)$
and compact subset $K $ of $\ov \Omega \smallsetminus \cO$, there holds
the inequality
\bea
\Vert u ; L_{\beta-1}^2(\Omega) \Vert + 
\Vert u ; L_\beta^2\partial (\Omega) \Vert \leq  
c\big( \Vert \nabla u ; L_{\beta}^2 (\Omega) \Vert +
\Vert u ; L^2(K ) \Vert \big) , \label{32}
\eea
where the constant $c$  depends on $\Omega$, $K$ and $\beta$ but not on $u$. 
\ENL

Proof. Replacing $u \mapsto r^{-\beta}u$ reduces the claim to the case $\beta=0$ 
which has been considered  in \cite[Sect.\,2]{na401}.  \ \ $\boxtimes$

\bigskip

We associate with the problem \eqref{31} the integral identity \cite{Lad}
\bea
(\nabla u, \nabla v)_\Omega - \lambda(u,v)_{\partial \Omega} = F(v) \hskip1cm
\forall \, v \in V_{-\beta}^1(\Omega) , \label{33}
\eea
 where $F \in V_{-\beta}^1(\Omega)^*$ is an (anti)linear functional
in $V_{-\beta}^1(\Omega)$, for instance,
\bea
F(v)= (f,v)_\Omega + (g,v)_{\partial \Omega} \ \ \mbox{with} \ 
f \in L_{\beta+1}^2 (\Omega) , \ g \in 
L_\beta^2(\partial\Omega). \label{33A}
\eea
According to \eqref{33A}, all terms in \eqref{33} are properly defined so 
that \eqref{33} defines a continuous mapping 
\bea
V_\beta^1(\Omega) \ni u \mapsto T_\beta(\lambda) u = F \in V_{-\beta}^1(\Omega)^*.
\label{34}
\eea
In \cite{na401} it is proven that the operator $T_0(\lambda)$ is Fredholm for 
$\lambda\in [0, \lambda_\dagger)$ but loses this property for $\lambda \geq 
\lambda_\dagger$. Notice that the latter fact follows from the failure of the 
inclusion $z^{\tau_\pm} \in V_0^1(\Pi^d)$ in the case \eqref{27}. We remark that $T_\beta(\lambda)$
is still Fredholm, if $\lambda\geq \lambda_\dagger$ and $\beta \not=0$, although 
this fact will be of no use here.

\subsection{Asymptotics of the solutions in the cusp.}
\label{sec2.4}
We consider problem \eqref{33} with the  right-hand side \eqref{33A},
where $\beta = -1$ and 
\bea
g=0, \ f \in L^2(\Omega).   \label{35}
\eea

The following assertion on the asymptotics of the solution of the Steklov
problem \eqref{31} can be found in \cite[Thm.\,2.6]{na549}. 

\BET
\label{T2}
Let $u \in V_1^1(\Omega)$ be a solution of the problem \eqref{34} with $\lambda \geq 
\lambda_\dagger$, $\beta =1$ and the right-hand side \eqref{33A}, \eqref{35}. Then,
$u$ has the representation
\bea
u(x) = \widetilde u (x) + \chi(x) \big( w(z) + W(z^{-2}y ,z) \big),
\label{36}
\eea
where the remainder $\widetilde u$ lives in $V_{-1}^1(\Omega)$ and  $\chi$ is a smooth 
cut-off function,
\bea
\chi = 1 \  \mbox{in} \ \Pi^{d/2}  \ and \ 
\chi = 0 \  \mbox{in} \ \Omega \smallsetminus \Pi^d  . \label{37}
\eea
Moreover, $w$ is the linear combination \eqref{23} with some coefficients $b_\pm$ 
depending on $f$ and including the functions \eqref{26} in the case $\lambda > 
\lambda_\dagger$ and \eqref{28A} in the case $\lambda= \lambda_\dagger$, while  
$W$ is determined by \eqref{22}, \eqref{23}, \eqref{24}, \eqref{29A}. 
Furthermore, there holds the  estimate 
\bea
\Vert \widetilde u ; V_{-1}^1(\Omega) \Vert + |b_+| + |b_-| \leq c 
\big( \Vert f ;L^2(\Omega) \Vert + \Vert u ; V_1^1(\Omega) \Vert \big)
\label{38}
\eea
with a coefficient $c$ independent of $f$ and $u$. 
\ENT

We emphasize that for $\lambda < \lambda_\dagger$, i.e., below the continuous spectrum
$\sigma_{\rm co}$, the asymptotic form of the solutions of the problem \eqref{31}
is different from that in Theorem \ref{T2}. 

\BER
\label{R2}
A direct calculation based on formulas \eqref{26}, \eqref{28A} and \eqref{22}
shows that the function $\chi(w +W)$ in Theorem \ref{T2} lives in $V_1^1(\Omega)$ 
but does not belong to $V_{-1}^1(\Omega)$, if $|b_+| + |b_-| \not=0$. Note that 
the solution $W$ of the Neumann problem \eqref{23}, \eqref{24} is determined up to the addentum
\bea
z^{- n + 3/2 \pm i \tau_0} W_\bullet \ \ \mbox{with} \ \tau_0 = \sqrt{ \Big(n -\frac32
\Big)^2 - \Lambda} \label{38A}
\eea
which is constant with respect to $\eta$. However, the term \eqref{38A} belongs 
to both spaces $V_{\pm 1}^1(\Pi^d)$ and thus can be omitted in \eqref{36}. This 
explains  why one requires the orthogonality condition \eqref{29A} on 
$W_0$. \ \ $\boxtimes$
\ENR

\BER
\label{R3}
A solution $u$ of the problem \eqref{33}  with data 
\eqref{35} belongs to the linear space $H_ {\rm loc}^2 (\ov \Omega \smallsetminus 
\cO)$. According to  \cite[Lem.\,3.2]{na549} [[correct reference?]], the second 
derivatives of the remainder  $\widetilde u \in V_{-1}^1(\Omega)$  belong to the 
space $L_1^2(\Omega)$ but not  to  $L_0^2(\Omega) = L^2(\Omega)$, although one 
might imagine so on the basis of \eqref{35}.
\ENR

\section{Radiation conditions and extension of the operator.}
\label{sec3}
\subsection{Generalized Green's formula.}
\label{sec3.1}
Given two right-hand sides $f^1$ and $f^2$  belonging to $L^2(\Omega)$,
let $u^1$ and $u^2$ be the  solutions of the problem \eqref{31}. Let 
also $b_\pm^1$  and  $b_\pm^2$ denote the coefficients in the linear combinations 
for $w^1$ and $w^2$ in \eqref{28}, which appear in the asymptotic formula 
\eqref{36} for $u^1$ and $u^2$,  respectively. We insert these solutions into the 
Green's formula on the truncated domain
$\Omega^\delta$, see \eqref{7A}. Passing to the limit $\delta \to +0$, we get
\bea
& & (f^1,u^2)_\Omega - u^1, f^2)_\Omega=  - \lim\limits_{\delta \to+0}
\big( (\Delta u^1, u^2)_{\Omega^\delta} -(\Delta u^2, u^1)_{\Omega^\delta} \big)
\roweq
\lim\limits_{\delta \to+0} \int\limits_{\omega^\delta} 
\big( \overline{u^2(y,\delta)} \partial_z  u^1(y,\delta) - 
u^1(y,s) \overline{\partial_z u^2(y,\delta) } \big) dy . \label{43}
\eea
First, we consider the case $\lambda > \lambda_\dagger$, when the entries of \eqref{28} are 
of the form \eqref{26}. We follow \cite[Sect.\,3.4]{na549}, see also
\cite{na493}, and use the decay properties of $\widetilde u^j (y,z)$ and 
$W^j(z^{-2} y,z )$, see \eqref{38}, \eqref{22}, to change in the limit
the integrand in \eqref{43} to 
\beas
\overline{w^2(\delta)} \partial_z w^1(\delta) - w^1(\delta) 
\overline{\partial_z w^2(\delta)}  .
\eeas
Hence, the representation \eqref{28}, \eqref{26}  of $w^0$ yields
\bea
& & (f^1, u^2)_\Omega - (u^1, f^2)_\Omega 
\roweq - w_0^2 \lim\limits_{\delta \to +0} \delta^{2(n-1)} |\omega|
\Big( \big( \overline{b_+^2 \delta^{\tau_+} + b_-^2 \delta^{\tau_-}} \big)
\big( \tau_+ b_+^1\delta^{\tau_+ -1 } + \tau_-b_-^1 \delta^{\tau_- -1} \big)
\rowmi
\big( b_+^1 \delta^{\tau_+} + b_-^1 \delta^{\tau_-} \big) 
\big( \overline{ \tau_+ b_+^2\delta^{\tau_+ -1 } + \tau_-b_-^2 \delta^{\tau_- -1} } \big)
\Big). \label{43A}
\eea
Thus, fixing the normalization factor in \eqref{26} as
\bea
w_0 = \frac{1}{\sqrt{2 |\omega|}} \Big( \Lambda - n+  \frac32\Big)^{-1/4} 
\hskip1cm  \mbox{for} \ \lambda > \lambda_\dagger, \label{44}
\eea
we derive the equality
\bea
- (\Delta u^1 , u^2) _\Omega + (u^1, \Delta u^2)_\Omega = 
i \overline{b_+^2} b_+^1 - i \overline{b_-^2} b_-^1 \label{45}
\eea
for functions of the form \eqref{37} satisfying the Steklov condition in \eqref{1}. 

The identity \eqref{45} holds true also in the case $\lambda = \lambda_\dagger$ 
with logarithmic singularities \eqref{28A}, when the normalization factor is chosen as
\bea
w_0 = \frac{1}{\sqrt{2 |\omega|}} .
\label{46}
\eea
This can be proven with calculations quite similar to \eqref{43}, \eqref{43A}. 

Taking into account Theorem \ref{T2} and generalizing the above calculations a 
bit (cf. \cite[Sect.\,3.4]{na549}) yields also the following assertion.

\BET
\label{T3}
Let $u^1$, $u^2 \in V_{-1}^1(\Omega) \cap H_{\rm loc}^2 (\overline{\Omega 
\smallsetminus \cO)}$  satisfy 
\bea
\Delta u^j \in L^2(\Omega) \ , \ \ \partial_\nu u^j - \lambda u^j \in 
L_{-1}^2(\partial \Omega), \ \ j =1,2. \label{Gj}
\eea
Then, these functions can be written in the form \eqref{36}, and 
there holds the  generalized Green's formula
\bea
& & - (\Delta u^1 , u^2) _\Omega + ( \partial_\nu u^1 - \lambda u^1, u^2)_{\partial \Omega}
+ (u^1, \Delta u^2)_\Omega
- 
(u^1 , \partial_\nu u^2 - \lambda u^2)_{\partial \Omega}
\nonumber \\ 
& & \ \ \ = i \overline{b_+^2} b_+^1 - i \overline{b_-^2} b_-^1 \label{Gg}
\eea
\ENT

\subsection{Wave processes in the cusp.}
\label{sec3.2}
We follow the paper \cite{na493}, which is related to a bit geometrically 
different  cuspidal irregularity of the boundary, and interpret the singular 
solutions \eqref{26} (detached in the right-hand side of \eqref{36}, 
$\lambda > 
\lambda_\dagger$) as waves  travelling along the cusp \eqref{0}. A clear 
physical  reason for such an interpretation can be found in the papers 
\cite{Mir, Kry, KoNa}  and others describing the Vibration Black Holes for 
acoustic and elastic waves.  The Mandelstam energy radiation principle  can be 
used to distinguish between  outgoing $w_+$ and incoming $w_-$ waves, namely, 
the former propagates to and the latter from the tip $\cO$; see \cite{Mand} 
and also \cite[Ch.\,5]{NaPl},  \cite{na493,KoNa}. As usual in scattering 
theory, this classification provides the following solution of the diffraction 
problem \eqref{1} in $\Omega$, see e.g. \cite{Wilcox,Mitta}, \cite[Ch.\,5]
{NaPl}, and Lemma \ref{lem3.1}, below:
\bea
Z(x) = \chi(x) \big( w_-(z) + W_-(\eta,z) \big) + s \chi(x)\big(w_+(z) + W_+(\eta,z) \big) 
+ \widetilde Z(x) \big) , \label{47}
\eea
which  is generated by the "incoming" wave $w_-$ and involves the 
scattering coefficient $s$ of the "outgoing" wave $w_+$.
The decomposition \eqref{47} is nothing but a concretization of \eqref{36}; the remainder
$ \widetilde Z$ belongs to 
$V_{-1}^1(\Omega)$ and $s$ is the so called scattering coefficient. Plugging the
harmonic function $Z$ into \eqref{45} gives
\bea
0 = i |s|^2 - i \ \ \Rightarrow \ \ s = e^{ i \Theta} \in \bbS^1 \subset \bbC.
\label{48}
\eea

Although we will provide a mathematical argument to support these formulas, it will
be convenient to use the physical terminology in the sequel. We will
write $Z(\lambda;x)$, $s(\lambda)$ and so on to indicate the dependence on the 
spectral parameter  $\lambda$.

\subsection{Weighted spaces with detached asymptotics.}
\label{sec3.3}
Let $\lambda \geq \lambda_\dagger $ and let $\cV^1(\Omega;\lambda)$ be the Banach space
composed of functions \eqref{36} and endowed with the norm
\bea
\Vert u ; \cV^1(\Omega ;\lambda) \Vert = 
\Vert \widetilde u ; V_{-1}^1(\Omega) \Vert + \sum_\pm |b_\pm| , \label{49}
\eea
where $\widetilde u$ is the remainder and $b_\pm$ are the coefficients of the 
linear combination \eqref{28}. Since $V_1^1(\Omega)^* \subset 
V_{-1}^1(\Omega)^*$, the operator
\bea
\cT(\lambda) : \cV^1 (\Omega; \lambda) \to V_1^1(\Omega)^* \label{50}
\eea
is nothing but the restriction of the operator $T_1(\lambda)$ to the 
subspace $\cV^1(\Omega;\lambda) \subset V_1^1(\Omega)$. In view of Theorem 
\ref{T2}, the operator \eqref{50} inherits the main properties of $T_1(\lambda)$, 
in  particular,  its kernel equals 
\beas
{\rm ker} \, \cT(\lambda) = {\rm ker} \, T_1(\lambda) =  \{
u \in V_1^1(\Omega) \, : \, T_1(\lambda) u=0 \}. 
\eeas
The operators $T_1(\lambda) $ and $T_{-1}(\lambda)$ are Fredholm and mutually adjoint,
and therefore
\bea
{\rm Ind}\, T_1(\lambda) = - {\rm Ind}\, T_{-1}(\lambda) . \label{51}
\eea
Furthermore, in view of Theorem \ref{T2}, their indices Ind\,$T_{\pm 1} (\lambda) = {\rm dim \, ker} \, T_{\pm 1}
(\lambda) - {\rm dim \, coker} \, T_{\pm 1}(\lambda)$ are related by
\bea
{\rm Ind}\, T_1(\lambda) = {\rm Ind}\, T_{-1}(\lambda) + 2 , \label{52}
\eea
where 2 is just the number of the free constants $b_\pm$ in the 
detached asymptotic term on the right-hand side of \eqref{36}. Obviously, 
ker\,$T_{-1}(\lambda) \subset {\rm ker} \, T_1(\lambda)$, hence, we can 
deduce from \eqref{51}, \eqref{52} that 
\bea
{\rm ker}\, T_1(\lambda) = \cZ \oplus {\rm ker}\, T_{-1}(\lambda) . \label{53}
\eea
where $\cZ$ is a subspace of dimension 1. 

\BEL
\label{lem3.1}
Let $\lambda\geq \lambda_\dagger$. The subspace $\cZ$ in \eqref{53} is spanned 
by the non-trivial solution $Z \in V_1^1(\lambda)$, see \eqref{47}, of the 
problem \eqref{33} with $F =0$, $\beta=1$. 
\ENL

Proof. A non-trivial element $Z \in \cZ$ has the form \eqref{36}, where $|b_+| +
|b_-| \not=0$ in the linear combination \eqref{28} (otherwise 
$Z  \in {\rm ker}\, T_{-1}(\lambda) \subset V_{-1}^1(\lambda)$). From \eqref{45} we 
deduce that $i |b_+|^2 - i|b_-|^2 =0$ so that none of the coefficients can vanish
and thus $Z$ indeed has the representation \eqref{47}.\ \ $\boxtimes$

\bigskip

The second component on the right of \eqref{53} consists of the so-called trapped
modes, i.e., solutions of the homogeneous Steklov problem 
\eqref{1} belonging to the space $V_{-1}^1(\Omega) \subset H^1(\Omega)$.
In \cite[Thm.\,2.6]{na549}  it was proven that ker\,$T_{-1}(\lambda) 
\subset  V_{-\beta}^1(\Omega)$ for any $\beta >0$, because the sum $w + W$ 
vanishes in  \eqref{36} and no other power-law terms appear. In other words,  the  
trapped  modes have at least superpower  decay rate as $x \to \cO$. 

Let $\cT_{\rm out}(\lambda)$ be the restriction of the operator \eqref{50}
to the subspace
\bea
\cV_{\rm out}^1(\Omega;\lambda) = \{ u \in \cV^1(\Omega;\lambda) \, :\, 
b_-=0\}. \label{54}
\eea
The condition on the right-hand side of \eqref{54} eliminates the incoming wave $w_-$
in the decomposition \eqref{36} so that $\cT_{\rm out}(\lambda)$ must be regarded
as the operator of the Steklov problem with the Mandelstam (energy) radiation 
conditions in the cusp (see, e.g., \cite[Ch.\,5]{NaPl}). 

Since Ind\,$\cT(\lambda) = 1$ by \eqref{51}, \eqref{52} and $Z \notin 
\cV_{\rm out}^1 (\Omega; \lambda)$, we observe that
\bea
{\rm ker} \, \cT_{\rm out} (\lambda) = {\rm ker}\, T_{-1}(\lambda) \subset
V_{-\beta}^1(\Omega)  \ \ \ \forall \, \beta \in \bbR   \label{55}
\eea
and that $\cT_{\rm out} (\lambda)$ is a Fredholm operator of index zero. 
Hence, problem \eqref{31} has a solution in the function space 
\eqref{54}, if and only if 
\beas
F(v) = 0 \ \ \forall \, v \in {\rm ker} \, \cT_{\rm out} (\lambda). 
\eeas
In this way, the Steklov problem with the Mandelstam radiation conditions
has all the general properties of traditional diffraction problems in 
cylindrical waveguides, cf. \cite{Mitta, Wilcox}.

\subsection{"Symmetric" realizations of the Steklov problem.}
\label{sec3.4}
We set, for $\theta \in [0, 2 \pi)$,
\bea
\cV_\theta^1(\Omega;\lambda) = \{ u \in \cV^1(\Omega;\lambda) \, : \, b_+ = e^{i \theta} 
b_-\}  \label{56}
\eea
and denote by $\cT_\theta^1(\lambda)$ the restriction of $\cT^1(\lambda)$ onto 
the  subspace \eqref{56}. Owing to \eqref{56}, formula \eqref{Gg} reads as 
\bea
& & - (\Delta u_\theta^1 , u_\theta^2) _\Omega + 
( \partial_\nu u_\theta^1 - \lambda u_\theta^1, u_\theta^2)_{\partial \Omega}
\roweq 
- (u_\theta^1, \Delta u_\theta^2)_\Omega
+ (u_\theta^1 , \partial_\nu u_\theta^2 - \lambda u_\theta^2)_{\partial \Omega}
\ \ \forall \, u_\theta^1, u_\theta^2 \in \cV_\theta^1(\Omega; \lambda).
 \label{57}
\eea
As this is the usual symmetric Green formula, we can regard $\cT_\theta(\lambda)$ as a 
symmetric Steklov problem operator, in contrast to the operator 
$\cT_{\rm out}(\lambda)$,  because for $u^1,u^2 \in \cV_{\rm out}^1(\Omega;\lambda)$ 
the right-hand side of \eqref{Gg} becomes $i \overline{b_+^2} b_+^1$, which does 
not vanish in general.  

Clearly,
\beas
{\rm ker} \, T_{-1}(\lambda) \subset {\rm ker} \,\cT_\theta(\lambda) 
\eeas
so that a trapped mode belongs to the kernel of $\cT_\theta(\lambda)$ 
for every parameter $\theta$. However, in the case
\beas
\theta =\Theta
\eeas
where $\Theta$ comes from \eqref{48} and \eqref{47}, ker\,$\cT_\theta(\lambda)$ coincides
with the subspace \eqref{53}, since the special solution \eqref{47}
with the scattering coefficient $s =e^{i\theta} = e^{i \Theta}$
belongs to the kernel of $\cT_\theta(\lambda)$. 

The above observations will be used in the next section to construct eigenvalues
belonging to \eqref{11}. In particular, the elements of ker\,$\cT_\theta(\lambda)$ will become, for some particular values of $\theta$, prototypes of the eigenfunctions of the singularly perturbed problem 
\eqref{8}--\eqref{10}.

\section{Spectrum in the domain with a blunted cusp.}
\label{sec4}
\subsection{Formal asymptotics.}
\label{sec4.1}
$1^\circ$. {\it Stable eigenvalues.}
We denote by $\lambda^{\sf tr}$ a number larger than $\lambda_\dagger$ and 
assume that there exists a non-zero element $u^{\sf tr}$ in 
ker\,$T_{-1}(\lambda^{\sf tr})$.  Since this trapped mode belongs to $V_{-\beta}^1(\Omega)$ 
for any $\beta>0$ and, therefore, leaves only a very small discrepancy in the Dirichlet 
condition  at the end $\omega^\varepsilon$ of $\Omega^\varepsilon$, the function $u^{\sf tr}$
is expected to be an excellent approximation of an eigenfunction of the problem 
\eqref{8}--\eqref{10}. Moreover, we will prove in Section \ref{sec4.3} that for 
some 
$\varepsilon (\lambda^{\sf tr})>0$ and all $\varepsilon \in (0, \varepsilon(\lambda^{\sf tr}))$,
there exists an eigenvalue $\lambda_{m(\varepsilon)}^\varepsilon$ in the sequence
\eqref{11} such that 
\bea
\big| \lambda_{m(\varepsilon)}^\varepsilon - \lambda^{\sf tr} \big| \leq c_\beta 
\varepsilon^\beta \ \ \ \forall\, \beta \in \bbR_+ , \label{61}
\eea
where $c_\beta$ is a constant independent of $\varepsilon$. In other words, the 
Steklov-Dirichlet problem in the domain \eqref{7A} with the blunted cusp
has an eigenvalue in the vicinity of the point $\lambda^{\sf tr} \in 
\sigma_{\rm co}$. 

Such a family of eigenvalues in $\Omega^\varepsilon$ with $\varepsilon \in (0, 
\varepsilon (\lambda^{\sf tr})]$ stays  close to a fixed point and has
the limit  $\lambda^{\sf tr}$ as $\varepsilon \to +0$ so that we call
them  "stable eigenvalues".

$2^\circ$. {\it Blinking eigenvalues.} Let us fix a point $\lambda^\flat > \lambda_\dagger$
and consider the solution $Z(\lambda^\flat ; \cdot)$ with the scattering 
coefficient 
\bea
s(\lambda^\flat) = e^{i \Theta (\lambda^\flat)} , \label{61A}
\eea
see \eqref{47} and \eqref{4}. According to \eqref{2}, the main asymptotic term
\bea
w_-(\lambda^\flat ; z) + s(\lambda^\flat ) w_+(\lambda^b;z) \label{61B}
\eea
in the decomposition of $Z(\lambda^\flat;z)$ vanishes at $z=\varepsilon$,
provided
\bea
& & \varepsilon^{-(n-3/2)- i\tau_0(\lambda^\flat)} + e^{i \Theta(\lambda^\flat)}
\varepsilon^{-(n-3/2)+ i\tau_0(\lambda^\flat)} = 0 
\row
\mbox{where} \ \ \tau_0(\lambda) = \sqrt{\frac{|\partial \omega|}{|
\omega|}\lambda - \Big( n - \frac32 \Big)^2}  , \label{62}
\eea
see \eqref{26}.
Thus,
\bea
- 2 \tau_0(\lambda^\flat) \ln \varepsilon = \Theta(\lambda^\flat) + \pi
\ \ \ (\mbox{mod} \ 2 \pi )   \label{62A}
\eea
and for the sequence $\{ \varepsilon_k^\flat  \}_{k=1}^\infty$, where
\bea
\varepsilon_k^\flat 
= e^{ -2(\tau_0(\lambda^\flat))^{-1} ( (2k+1)\pi + \Theta(\lambda^\flat))}  
\to 0 \ \ \mbox{as} \ k \to + \infty , \label{63}
\eea
the discrepancy left by the function $Z(\lambda^\flat ; \cdot)$ to \eqref{10} is 
small.  Furthermore, we will prove that, for all large $k$, the problem 
\eqref{8}--\eqref{10} in $\Omega^{\varepsilon_k^\flat}$ has an eigenvalue $
\lambda_{m_k}^{\varepsilon_k^\flat}$ such that
\bea
\big| \lambda_{m_k}^{\varepsilon_k^\flat} - \lambda^\flat \big| \leq c_\flat 
\varepsilon_k^\flat |\ln \varepsilon_k^\flat |^{-1/2}  , \label{64}
\eea
where $c_\flat$ is independent of $k$. 

We call $\lambda^\flat$ a "blinking eigenvalue" for the following reason:
when $\varepsilon \to +0$, there emerges an eigenvalue of the problem 
\eqref{8}--\eqref{10} in the vicinity of the point $\lambda^\flat$ for values 
of  $\varepsilon$ obeying the period $\pi \tau_0(\lambda^\flat)^{-1}$ in the 
logarithmic scale $|\ln \varepsilon|$. By the argument at  the end of 
$3^\circ$, below, the point $\lambda^\flat > \lambda_\dagger$ itself  becomes
a true eigenvalue of the  problem \eqref{8}--\eqref{10} for some
$\varepsilon$ close to any $\varepsilon_k^\flat$ of \eqref{63}: 
we emphasize that  {\it every} point $\lambda^\flat > \lambda_\dagger$ 
becomes such a blinking eigenvalue. This observation also allows us to 
construct in  Section \ref{sec5.1} a singular Weyl sequence for the operator 
$\cS$, \eqref{4}, at any point $\mu \in (0, \mu_\dagger)$.

$3^\circ$. {\it Gliding eigenvalues.} Since the eigenvalues of the problem 
\eqref{8}--\eqref{10} depend continuously on the small parameter $\varepsilon 
>0$, see e.g. \cite[Ch.\,7,\,Sec.\,6.5]{Kato}, the effect of "blinking" 
ought to cause them move fast along the semi-axis $(\lambda_\dagger, + \infty)$ 
as functions of $\varepsilon$, cf. the papers \cite{na564} and \cite{na595},
which deal with spectral problems for differential operators with sign-changing
coefficients. We give a hypothetical explanation of this phenomenon on the level 
of formal asymptotic analysis. To this end, we compute the derivative 
$\partial_\varepsilon \lambda^\flat (\varepsilon_k^\flat)$ from the equation 
\eqref{62A}, and using \eqref{62}, obtain
\bea
\frac{\partial \lambda^\flat}{\partial \varepsilon} (\varepsilon_k^\flat) =
\frac{2}{\varepsilon_k^\flat} 
\big( \lambda^\flat  (\varepsilon_k^\flat) -\lambda_\dagger \big)  
\Big( \frac1{ |\ln \varepsilon_k^\flat|}
+ O\Big( \frac1{ |\ln \varepsilon_k^\flat|^2} \Big) \Big) .\label{epla}
\eea
Formula \eqref{epla} demonstrates the rapid "fall" at a distance from
$\lambda_\dagger$ and the smooth "landing" of it at the 
threshold. (The eigenvalues could be described
as parachutists releasing their chutes only very close to the 
surface.)

Furthermore, since the eigenvalues depend continuously on the parameter
$\varepsilon$, we observe that the gliding eigenvalues descending along the 
interval $(\lambda_\dagger, + \infty)$ must cross every point $\lambda
> \lambda_\dagger$. Hence, every $\lambda \in (\lambda_\dagger, + \infty)$
becomes  a true eigenvalue of the Steklov-Dirichlet problem 
\eqref{8}--\eqref{10} for some $\varepsilon$. By the argument in 
$2^\circ$, this happens almost periodically in the $|\ln \varepsilon|$-scale, 
that is, for  $\varepsilon $ very close to the computed values \eqref{63}.

$4^\circ$. {\it The threshold case.} According to \eqref{epla}, the
threshold $\lambda_\dagger$ absorbes all "gliding eigenvalues" in the limit
$\varepsilon \to +0$. However, there is no "blinking" phenomenon
related with $\lambda_\dagger$. 
Indeed, according to \eqref{28A} the equality $Z(\lambda_\dagger; \varepsilon,y) 
=0$ yields 
\beas
w_0\varepsilon^{-n + 3/2} (1 - i \ln \varepsilon) + e^{i \Theta(\lambda_\dagger)}
w_0\varepsilon^{-n + 3/2} (1 +i \ln \varepsilon) = 0
\eeas
and hence
\bea
e^{i \Theta(\lambda_\dagger) } = - \frac{1- i \ln \varepsilon}{1 + i \ln \varepsilon}.
\label{thep}
\eea
Since the right-hand side of \eqref{thep} tends to 1  as $\varepsilon \to +0$,
the blinking eigenvalues do not occur at all, and, moreover, a near-threshold
eigenvalue may appear in the special situation $\Theta(\lambda_\dagger) =0$
only.

\subsection{Operator formulation of the problem in $\Omega^\varepsilon$.}
\label{sec4.2}
We define the Hilbert space $\cH^\varepsilon$, which consists of functions 
$u^\varepsilon \in H^1(\Omega^\varepsilon)$ satisfying the Dirichlet condition 
\eqref{10}, and endow it with the scalar product 
\bea
\langle u^\varepsilon , v^\varepsilon \rangle_\varepsilon =
(\nabla u^\varepsilon , \nabla v^\varepsilon )_{\Omega^\varepsilon} + 
( u^\varepsilon ,  v^\varepsilon )_{\partial \Omega^\varepsilon} .\label{F1}
\eea
The operator $\cS^\varepsilon \cH^\varepsilon \to \cH^\varepsilon$, defined by
\bea
\langle \cS^\varepsilon u^\varepsilon , v^\varepsilon \rangle_\varepsilon
= ( u^\varepsilon ,  v^\varepsilon )_{\partial \Omega^\varepsilon} 
\ \ \forall \, u^\varepsilon, v^\varepsilon \in \cH^\varepsilon \label{F2}
\eea
is positive, symmetric, and continuous, therefore self-adjoint. In view of 
\eqref{F1} and \eqref{F2}, the variational formulation of the Steklov-Dirichlet problem
\eqref{8}--\eqref{10} reads as 
\bea
(\nabla u^\varepsilon , \nabla v^\varepsilon )_{\Omega^\varepsilon} =
\lambda^\varepsilon( u^\varepsilon ,  v^\varepsilon )_{\partial \Omega^\varepsilon} \ \ 
\forall \,  v^\varepsilon \in \cH^\varepsilon ,
\label{F3}
\eea
and it converts to the abstract equation 
\bea
\cS^\varepsilon u^\varepsilon = \mu^\varepsilon u^\varepsilon \ \ 
\mbox {in} \ \cH^\varepsilon  , \label{F4}
\eea
where the spectral parameters are related in the same way as in \eqref{5}. The 
surface $\partial \Omega^\varepsilon$ is Lipschitz and thus the operator 
$\cS^\varepsilon$ is compact, hence, as well known,  
the essential spectrum of $\cS^\varepsilon$ consists of the single point 
$\mu^\varepsilon =0$ and the discrete spectrum of the sequence 
$\{\mu_p^\varepsilon \}_{p \in \bbN} \subset (0,1)$ convergent to 0. The 
sequence turns into \eqref{11} via the formula
\eqref{5}.

The next assertion is known as the lemma on "near eigenvalues", cf. \cite{ViLu},
and it is a direct consequence of the spectral decomposition of the resolvent, see
e.g. \cite[\S\,6.2.]{BiSo}. 

\BEL
\label{LNE}
Let $U^\varepsilon \in \cH^\varepsilon$ and $M^\varepsilon > 0$ be such that
\bea
\Vert U^\varepsilon ; \cH^\varepsilon \Vert = 1 \ \ \mbox{and} \ \
\Vert \cS^\varepsilon U^\varepsilon - M^\varepsilon 
U^\varepsilon; \cH^\varepsilon\Vert
=: \delta^\varepsilon \in (0,M^\varepsilon). \label{F5}
\eea
Then, the operator $S^\varepsilon$ has an eigenvalue $\mu_p^\varepsilon$ such that
\bea
|M^\varepsilon- \mu_p^\varepsilon| \leq \delta^\varepsilon. \label{F6}
\eea
\ENL

It will be important in the sequel that if  the condition 
\bea
\frac{\delta^\varepsilon}{M^\varepsilon} \leq \frac12 \label{F7}
\eea
holds, then the relations \eqref{F6} and \eqref{5} imply
\bea
& & \Big| 1 + \lambda_p^\varepsilon- \frac{1}{M^\varepsilon}\Big|
\leq \frac{\delta^\varepsilon}{M^\varepsilon}  (1 + \lambda_p^\varepsilon) 
\nonumber \\
\Rightarrow & & 
1 + \lambda_p^\varepsilon \leq \frac{2}{M^\varepsilon}
\ \ \Rightarrow \ \ 
\Big| 1 + \lambda_p^\varepsilon - \frac{1}{M^\varepsilon} \Big| \leq 
\frac{2 \delta^\varepsilon}{(M^\varepsilon)^2}. \label{F8}
\eea

\subsection{Justification of the "stable asymptotics".}
\label{sec4.3}
We assume that $u^{\sf tr} \in {\rm ker} \, T_{-1}(\lambda^{\sf tr}) \smallsetminus \{ 0 \}$ for
some $\lambda^{\sf tr} > \lambda_ \dagger$, cf. Section \ref{sec4.1},$1^\circ$, and set
\bea
M= (1 + \lambda^{\sf tr})^{-1} \ , \ \ U^\varepsilon = \Vert X^\varepsilon u^{\sf tr} ; 
\cH^\varepsilon \Vert^{-1} X^\varepsilon u^{\sf tr} , \label{F9}
\eea
where $X^\varepsilon$ is the smooth cut-off function 
\bea
& & X^\varepsilon(x) = 1, \,  x \in \Omega^\varepsilon  
\smallsetminus \Pi^{3 \varepsilon}, \ \ \ \ 
X^\varepsilon(x) = 0, \, x \in \Pi^{2 \varepsilon} ,
\row
\big| \nabla X^\varepsilon (x) \big| \leq c\varepsilon^{-1}. \label{F10}
\eea
Since a non-zero harmonic function cannot vanish on a set of positive $n$-measure
and $\lambda^{\sf tr}$ is positive, we have
\beas
\Vert X^\varepsilon u^{\sf tr} ; \cH^\varepsilon \Vert 
\geq\big( \Vert \nabla u^{\sf tr} ; L^2(\Omega \smallsetminus \Pi^d) \Vert^2 +
\Vert  u^{\sf tr} ; L^2(\partial \Omega \smallsetminus \Gamma^d) \Vert^2 \big)^{1/2}
\geq c_u > 0.
\eeas

Let us evaluate the quantity $\delta^\varepsilon$ in \eqref{F5}. Using the definition
of the Hilbert norm, we apply \eqref{F1}, \eqref{F2}, \eqref{F9} and write
\bea
& & \delta^\varepsilon = \sup \big| \langle \cS^\varepsilon U^\varepsilon - M^\varepsilon 
U^\varepsilon, V^\varepsilon \rangle_\varepsilon \big|
\roweq
(1 + \lambda^{\sf tr})^{-1} \Vert X^\varepsilon u^{\sf tr} ; \cH^\varepsilon\Vert^{-1}
\sup \big|\big( \nabla(X^\varepsilon u^{\sf tr} ),\nabla V^\varepsilon 
\big)_{\Omega^\varepsilon} - \lambda^{\sf tr} (X^\varepsilon u^{\sf tr}, V^\varepsilon
)_{\partial \Omega^\varepsilon} \big|  \label{F11}
\eea
where the supremum is taken over the unit ball of $\cH^\varepsilon$, i.e.,
$\Vert V^\varepsilon; \cH^\varepsilon \Vert \leq 1$. 

We have
\bea
& & \big( \nabla(X^\varepsilon u^{\sf tr}), \nabla V^\varepsilon \big)_{\Omega^\varepsilon}
- \lambda^{\sf tr} \big( X^\varepsilon u^{\sf tr},  V^\varepsilon \big)_{
\partial \Omega^\varepsilon} 
\roweq
\big( \nabla u^{\sf tr}, \nabla(X^\varepsilon V^\varepsilon ) \big)_{\Omega}
- \lambda^{\sf tr} \big(  u^{\sf tr},  X^\varepsilon V^\varepsilon \big)_{
\partial \Omega} 
\rowpl
\big(  u^{\sf tr} \nabla X^\varepsilon, \nabla  V^\varepsilon \big)_{
\Pi^{3 \varepsilon} \smallsetminus \Pi^{2 \varepsilon} } 
- \big( \nabla u^{\sf tr} , V^\varepsilon \nabla  X^\varepsilon \big)_{
\Pi^{3 \varepsilon} \smallsetminus \Pi^{2 \varepsilon}} . \label{F12}
\eea
Here, the changes of the integration domains $\Omega^\varepsilon$, $\partial 
\Omega^\varepsilon$ and $\Omega^\varepsilon$, respectively, to
$\Omega$, $\partial \Omega$ and $\Pi^{3 \varepsilon} \smallsetminus \Pi^{2 
\varepsilon}$ can be made because of the properties \eqref{F10} of the 
cut-off function  $X^\varepsilon$. In particular, the product $X^\varepsilon 
V^\varepsilon$ falls into
the space $\cH$ with the norm \eqref{2} so that the sum of the two terms
on the right-hand side of \eqref{F12} vanishes due to integral identity
\eqref{3}, where we put $\lambda= \lambda^{\sf tr}$ and $u= u^{\sf tr}$.
Since $u^{\sf tr} \in V_ {-\beta}^1(\Omega)$  for any $\beta$, Lemma \ref{L1}
leads to the estimates
\bea
& & \big| ( u^{\sf tr} \nabla X^\varepsilon, \nabla V^\varepsilon)_{
\Pi^{3 \varepsilon} \smallsetminus \Pi^{2 \varepsilon}} \big| \leq 
c \varepsilon \Vert u^{\sf tr} ; L^2(\Pi^{3 \varepsilon} \smallsetminus 
\Pi^{2 \varepsilon} ) \Vert \, \Vert \nabla V^\varepsilon ; 
L^2(\Omega^\varepsilon) \Vert
\rowleq
c_\beta \varepsilon^{-1} \varepsilon^{\beta+1} \Vert u^{\sf tr} ; L_{-\beta-1} 
(\Omega) \Vert
\, \Vert V^\varepsilon; \cH^\varepsilon\Vert \leq c_\beta \varepsilon^\beta, 
\row
\big| ( \nabla u^{\sf tr} , V^\varepsilon\nabla X^\varepsilon)_{
\Pi^{3 \varepsilon} \smallsetminus \Pi^{2 \varepsilon}} \big| 
\rowleq 
c\Vert \nabla u^{\sf tr} ; 
L^2(\Pi^{3 \varepsilon} \smallsetminus \Pi^{2 \varepsilon} ) \Vert \, 
\varepsilon^{-1} \Vert V^\varepsilon ; 
L^2(\Pi^{3 \varepsilon} \smallsetminus \Pi^{2 \varepsilon}) \Vert
\rowleq
c_\beta\varepsilon^{\beta}  \Vert \nabla  u^{\sf tr} ; L_{-\beta} (\Omega) \Vert
\, \Vert V^\varepsilon; L_{-1}^2 (\Omega^\varepsilon) \Vert 
\leq 
c_\beta \varepsilon^\beta \Vert V^\varepsilon; \cH^\varepsilon\Vert 
\leq c_\beta \varepsilon^\beta . \label{F13}
\eea
Thus, $\delta^\varepsilon$ does not exceed $c_\beta\varepsilon^\beta$ for any $\beta \in 
\bbR$, and we have proven the following assertion.

\BET
\label{ASY1}
Assume that $\lambda^{\sf tr} > \lambda_\dagger$ and that there exists a trapped mode 
$u^{\sf tr} \not= 0$ belonging  to {\rm ker}\,$T_{-1}(\lambda^{\sf tr}) $.
Then, for some $\varepsilon^{\sf tr} >0$, the formula \eqref{61} is valid for
all $\varepsilon \in (0, \varepsilon^{\sf tr}]$ and for some eigenvalue 
$\lambda_{m(\varepsilon)}^\varepsilon$ belonging to the spectrum \eqref{11}
of the Steklov-Dirichlet problem \eqref{8}--\eqref{10}.
\ENT

\subsection{Boundary layer.}
\label{sec4.4}
The asymptotic structure of an approximate eigenfunction, which is based on the
function $Z(\lambda^\flat; x)$, $\lambda^\flat > \lambda_\dagger$, is much
more complicated than that in the previous section. Although the main term
\eqref{61B} of the expansion \eqref{47} vanishes at the end 
$\omega^{\varepsilon_k^\flat}$ of the blunted domain $\Omega^{\varepsilon_k^\flat}$ 
(with the small parameter \eqref{63}) and the remainder $Z(\lambda^\flat; \cdot)
\in V_{-1}^1(\Omega)$ decays sufficiently fast, the correction term
\bea
W_-(\lambda^\flat ;z) + s(\lambda^\flat) W_+(\lambda^\flat ; z) \label{Y1}
\eea
leaves a significant discrepancy in the Dirichlet condition \eqref{10}. To 
compensate this we follow a general approach of \cite[Ch.\,15]{MaNaPl} and 
employ  the stretched coordinates
\bea
x=(y,z) \mapsto \xi = (\xi', \xi_n) = \big( \varepsilon^{-2} y ,
\varepsilon^{-2} (z- \varepsilon) \big) \label{Y2}
\eea
to describe the boundary layer phenomenon. We emphasize the apparent difference 
between  the transversal coordinates in \eqref{Y2} and \eqref{0}. By \eqref{0},
\eqref{7}, using  \eqref{Y2} and setting formally $\varepsilon =0$
convert the domain \eqref{7A} to the half-cylinder 
\bea
\Xi = \omega \times \bbR_+. \label{Y3}
\eea
In this transformation the Laplacian gets the factor $\varepsilon^{-4}$ and,
by \eqref{23A}, the  Steklov condition \eqref{9} reduces asymptotically
to  the Neumann one  on the lateral side $\Sigma = \partial \omega \times \bbR_+$ of 
\eqref{Y3}. Consequently, the problem for the boundary layer reads as
\bea
& & - \Delta_\xi Y_{\pm} (\xi) = 0.  \,   \xi \in \Xi,  \ \ \ \ 
\partial_{\nu'}  Y_{\pm} (\xi) = 0, \,  \xi \in \Sigma,
\row \ \ \, 
Y_\pm(\xi',0) = W_\pm(\xi',1), \, \xi' \in \omega. \label{Y4}
\eea
The Fourier method proves that, under the orthogonality condition \eqref{29A},
problem \eqref{Y4} has a unique solution with exponential decay at infinity,
namely
\bea
e^{\beta_\Xi \xi_n} Y_\pm \in H^2(\Xi) \ \ \mbox{for some} \ \beta_\Xi > 0. 
\label{Y4a}
\eea
Notice that also the second derivatives of $Y_\pm$ belong to $L^2(\Xi)$, since
there are no "strong" singularities at the Dirichlet-Neumann collision corner
point of opening $\pi/2$, see, e.g., \cite[Ch.\,2\,and\,11]{NaPl}.

\subsection{Justification of the "blinking asymptotics".}
\label{sec4.5}
Fixing some $\lambda^\flat > \lambda_\dagger$, we consider Steklov-Dirichlet problem
\eqref{8}--\eqref{10} in the domain $\Omega^{\varepsilon_k^\flat}$ with the
small $\varepsilon_k^\flat$ as in \eqref{63}, and the  solution 
$Z(\lambda^\flat;x)$ with the scattering coefficient \eqref{61A}.

Using \eqref{5} and \eqref{47} we write 
\bea
& & M^\varepsilon 
=(1+ \lambda^\flat)^{-1} , \ U^\varepsilon= \Vert u^\flat ; \cH^\varepsilon 
\Vert^{-1} u^\flat, 
\row
u^\flat (x) = X^\varepsilon(x) \widetilde Z(\lambda^\flat; x) 
+ 
\chi(x) \sum_\pm s_\pm(\lambda^\flat) \big( w_\pm(z) \nonumber \\
& & \hskip1.5cm 
+ W_\pm(z^{-2} y,z) -
\varepsilon^{2-(n-3/2) \pm i \tau_0(\lambda^\flat)}
Y_\pm(z^{-2} y ,\varepsilon^{-1} z)\big) , \label{Y5}
\eea
where we set  $s_-(\lambda^\flat ) =1$ and $s_+(\lambda^\flat) =
 s(\lambda^\flat)$ to shorten the notation. 

First of all, we evaluate the norm
\beas
\Vert u^\flat ; \cH^\varepsilon \Vert \geq 
\Vert u^\flat ; L^2(\partial \Omega^\varepsilon \smallsetminus 
\omega
\Vert. 
\eeas
Owing to the basic properties of $\widetilde Z, W_\pm$ and $Y_\pm$,
we have
\beas
\Vert u^\varepsilon ; L^2(\partial \Omega^\varepsilon \smallsetminus 
\omega
\Vert \geq J^\varepsilon - c ,
\eeas
where 
\bea
J^\varepsilon &=&
\int\limits_{\varepsilon_k^\flat}^d \int\limits_{\partial\omega^z}
\big| w_-(\lambda^\flat; z) + e^{i \Theta(\lambda^\flat)}
w_+(\lambda^\flat;z) \big|^2 ds_y dz
\roweq
\int\limits_{\varepsilon_k^\flat}^d
z^{2(n-2)} |\partial \omega| z^{-2n +3}
\big| 1+  e^{i ( \Theta(\lambda^\flat) + 2 \tau_0(\lambda^\flat))z}
 \big|^2  dz
\roweq
2 |\partial \omega|
\int\limits_{\varepsilon_k^\flat}^d 
\big( 1+ \cos\big(  ( \Theta(\lambda^\flat) + 2 \tau_0(\lambda^\flat) z\big)
\big) \frac{dz}{z}  
\roweq
2 |\partial \omega|
\big( \ln z + {\rm Ci} \big(  ( \Theta(\lambda^\flat) + 2 \tau_0(\lambda^\flat) 
z\big) \big) \Big|_{\varepsilon_k^\flat}^d \geq
 c_k^\flat \big|\ln \varepsilon_k^\flat \big| ,
\eea
where 
\beas
{\rm Ci}(\tau) = - \int\limits_\tau^{\infty} \frac{\cos t}{t} dt.
\eeas
is the cosine integral function. Thus,
\bea
\Vert u^\flat ; \cH^\varepsilon \Vert \geq C_k^\flat 
\big|\ln \varepsilon_k^\flat \big|^{1/2}\ \ , \ \ C_k^\flat > 0.
\label{big}
\eea

By \eqref{F10}, \eqref{62} and \eqref{Y4} we see that $u^\varepsilon=0$ at
$z= \varepsilon$. We continue the calculation of the quantity 
$\delta^\varepsilon$ in \eqref{F5} as follows:
\bea
& & \delta^\varepsilon = ( 1 + \lambda^\flat)^{-1} 
\Vert u^\flat ; \cH^\varepsilon \Vert^{-1} \sup \big|
(\nabla u^\flat, \nabla V^\varepsilon)_{ \Omega^\varepsilon} -
\lambda^\flat ( u^\flat,  V^\varepsilon)_{\partial \Omega^\varepsilon}  \big|
\roweq
( 1 + \lambda^\flat)^{-1} 
\Vert u^\flat ; \cH^\varepsilon \Vert^{-1} \sup \big|
(\Delta u^\flat, V^\varepsilon)_{ \Omega^\varepsilon} 
- (\partial_\nu u^\flat - \lambda^\flat u^\flat ,  V^\varepsilon)_{ 
\partial \Omega^\varepsilon} \big| \label{Y6}
\eea
where again the supremum is taken over the unit ball of $\cH^\varepsilon$.
We denote by $I(V^\varepsilon)$ the expression inside the last moduli of 
\eqref{Y6} and write it as the sum $I_Z(V^\varepsilon) + I_Y(V^\varepsilon) +
I_X(V^\varepsilon) $, where
\bea
 I_Z(V^\varepsilon) & = & (\Delta Z(\lambda^\flat; \cdot), 
 X^\varepsilon V^\varepsilon)_{ \Omega} 
- \big( (\partial_\nu - \lambda^\flat)  Z(\lambda^\flat; \cdot), 
 X^\varepsilon V^\varepsilon \big)_{ 
\partial \Omega^\varepsilon} ,
\nonumber \\
I_Y(V^\varepsilon)& =& \sum_\pm s_\pm(\lambda^\flat) \varepsilon^{2- (n-3/2) 
\pm i \tau_0(\lambda^\flat)} 
\Big( (\Delta(\chi Y_\pm), V^\varepsilon)_{\Omega^\varepsilon}  
- \big( (\partial_\nu -\lambda^\flat) \chi Y_\pm  , 
V^\varepsilon\big)_{\partial \Omega^\varepsilon} \Big) ,
\nonumber \\
I_X(V^\varepsilon) &=&  \Big(\Delta \big((1- X^\varepsilon) \widetilde 
Z(\lambda^\flat ; \cdot) \big) , V^\varepsilon\Big)_{ \Omega^\varepsilon}
- \big( (\partial_\nu - \lambda^\flat )(1-X^\varepsilon)
\widetilde  Z(\lambda^\flat ; \cdot), V^\varepsilon\big)_{\partial 
\Omega^\varepsilon}
\roweq
\lambda^\flat \big( (1- X^\varepsilon) \widetilde 
Z(\lambda^\flat ; \cdot)  , V^\varepsilon \big)_{\partial \Omega^\varepsilon}
- \big( \nabla ( (1-X^\varepsilon)
\widetilde  Z(\lambda^\flat ; \cdot)), \nabla 
V^\varepsilon\big)_{\Omega^\varepsilon} .  \label{Y7}
\eea
Clearly, $I_Z(V^\varepsilon) =0$ since \eqref{47} is a  
solution of \eqref{1}. Moreover, using Lemma \ref{L1} we obtain, 
similarly to \eqref{F13} with $\beta=1$,
\bea
|I_X(V^\varepsilon)| &\leq &
c\Big( \lambda^\flat  \Vert \widetilde Z  ; L^2(\Gamma^{3 \varepsilon} 
\smallsetminus  \Gamma^{2 \varepsilon} ) \Vert \, 
\Vert V^\varepsilon ; L^2(\partial \Omega^\varepsilon) \Vert
\rowpl
\big( \Vert \nabla \widetilde Z  ; L^2(\Pi^{3 \varepsilon} 
\smallsetminus  \Pi ^{2 \varepsilon} ) \Vert 
+ \varepsilon^{-1}  \Vert \widetilde Z  ; L^2(\Pi^{3 \varepsilon} 
\smallsetminus  \Pi ^{2 \varepsilon} ) \Vert \big) 
\Vert \nabla V^\varepsilon ; L^2(\partial \Omega^\varepsilon) \Vert
\Big)
\rowleq
c\varepsilon \big( \Vert \widetilde Z ; L_{-1}^2 (\partial \Omega) \Vert
+\Vert \widetilde Z ; V_{-1}^1 ( \Omega) \Vert
\big) \, \Vert V^\varepsilon; \cH^\varepsilon\Vert \leq c \varepsilon. 
\label{Y8}
\eea
To estimate the remaining term $I_Y(V^\varepsilon)$, we
write
\beas
& & J_\pm (x)  = \Delta\big( \chi(z) Y_\pm(\eta,\xi_n)\big) -
[\Delta, \chi(z) ]Y_\pm(\eta, \xi_n)
\roweq
\chi(z) z^{-4} \Delta_\eta Y_\pm (\eta, \xi_n) + \chi(z) 
\big( \varepsilon^{-2} \partial_{\xi_n} - 2 \eta z^{-1} \nabla_\eta
\big)^2Y_\pm (\eta,\xi_n). 
\eeas
Owing to the definition \eqref{37} of $\chi$, the supports of the coefficient
functions in the commutator $[\Delta, \chi]$ are included in the set
$\{ x =(y,z) \in \overline{ \Pi^d} \, : \, z \geq d/2 \}$,
where the functions $Y_\pm$ are exponentially small. Hence,
\beas
\big| \varepsilon^{2 - (n -3/2) \pm i \tau_0(\lambda^\flat) }
\big( [\Delta, \chi] Y_\pm , V^\varepsilon\big)_{\Pi^d } \big| 
\leq Ce^{-\delta / \varepsilon^2}  , \ \delta > 0.
\eeas
Furthermore, 
\beas
|J_\pm(x) | &\leq& c \Big( (\varepsilon^{-4} - z^{-4} ) |\Delta_\eta
Y_\pm(\eta, \xi_n)| 
\rowpl
z^{-1}(\varepsilon^{-2} + z^{-1} ) 
\big| \nabla_{(\eta, \xi_n)}^2 Y_\pm (\eta, \xi_n)\big| +
z^{-2} |\nabla_\eta Y_\pm (\eta,\xi_n)| \Big) 
\eeas
because $Y_\pm$ is harmonic. Thus, the evident relations
\beas
z^{T_1} e^{-\beta_\Xi \varepsilon^{-2} (z- \varepsilon)}
\leq c_1 \varepsilon^{T_1}\ , \ 
(z-\varepsilon)^{T_2} e^{-\beta_\Xi \varepsilon^{-2} (z- \varepsilon)}
\leq c_2 \varepsilon^{2 T_1}
\eeas
and $dx = z^{2(n+1) } \varepsilon^2 d\eta d\xi_n$ imply that 
\beas
& & \big| \varepsilon^{2 - (n-3/2)\pm i\tau_0(\lambda^\flat)}
(J_\pm , V^\varepsilon)_{\Pi^d \smallsetminus \Pi^\varepsilon} \big|
\rowleq
c \varepsilon^{-n +7/2} \Big( \int\limits_{\Pi^d \smallsetminus \Pi^\varepsilon}
z^2 |J_\pm(x)|^2 dx \Big)^{1/2}
\Vert z^{-1} V^\varepsilon; L^2(\Pi^d \smallsetminus \Pi^\varepsilon) \Vert
\rowleq
c \varepsilon^{-n +7/2} \Big( \max\limits_{ \varepsilon \leq z \leq d}
\big( (\varepsilon^{-8} (z- \varepsilon)^2 + \varepsilon^{-4} + z^{-4} )
z^{2(n-1)} e^{-2 \beta_\Xi \varepsilon^{-2} (z-\varepsilon)} \big) 
\nonumber \\
& & \times
\varepsilon^2 \int\limits_\Xi e^{2 \beta_\Xi \xi_n} \big( |\nabla_\xi^2
Y_\pm (\xi) |^2 + |\nabla_\xi Y_\pm(\xi)|^2 \big) d\xi \Big)^{1/2}
\Vert V^\varepsilon; \cH^\varepsilon\Vert 
\leq 
c \varepsilon^{3/2} .
\eeas
The surface integrals in $I_Y(V^\varepsilon)$ are treated in a similar way.
Using \eqref{23A} and \eqref{Y4} we obtain
\beas
& & \big| \varepsilon^{2 - (n-3/2)\pm i\tau_0(\lambda^\flat)}
\big( (\partial_\nu - \lambda^\flat) \chi Y_\pm , V^\varepsilon)_{\Gamma^d 
\smallsetminus \Gamma^\varepsilon} \big|
\rowleq
c \varepsilon^{-n +7/2} \Big( e^{-\delta \varepsilon^{-2} }+ 
\max\limits_{ \varepsilon \leq z \leq d}
\Big( (1+ \varepsilon^{-4} z^2 ) e^{-2 \beta_\Xi \varepsilon^{-2} (z-\varepsilon)} z^{2(n-2)} \Big) 
\nonumber \\
& & \times
\varepsilon^2 \int\limits_\Sigma e^{2 \beta_\Xi \xi_n} \big( |\nabla_\xi
Y_\pm (\xi) |^2 + | Y_\pm(\xi)|^2 \big) d\xi \Big)^{1/2}
\Vert V^\varepsilon; L^2(\Gamma^d \smallsetminus  \Gamma^\varepsilon )
\Vert 
\leq c \varepsilon^{3/2} .
\eeas
Collecting all the estimates and using Lemma \ref{LNE} yields the
following, desired assertion. 

\BET
\label{ASY2}
Let $\lambda^\flat > \lambda_\dagger$ and $\varepsilon_k^\flat$ be the small
parameter \eqref{63} such that \eqref{62A} holds for the scattering coefficient 
$e^{i \Theta(\lambda^\flat)}$ in the solution \eqref{47} of the problem \eqref{1}
in $\Omega$. Then, the problem \eqref{8}--\eqref{10} in $\Omega^\varepsilon$
has an eigenvalue $\lambda_{m_k}^{\varepsilon_k^\flat}$ satisfying the 
inequality \eqref{64}. 
\ENT

\section{Concluding remarks.}
\label{sec5}
\subsection{Singular Weyl sequence.}
\label{sec5.1}
Let the sequence   $ \{ \varepsilon_k^\flat \}_{k=1}^\infty$ be as in
\eqref{63}. 
We define the functions $u_k^\flat$ by using formula \eqref{Y5}
and extend them as zero from $\Omega^{\varepsilon_k^\flat}$ 
to the entire domain  $\Omega$. Let us show that the functions
\bea
U_k^\flat = \Vert u_k^\flat ; \cH\Vert^{-1} 
u_k^\flat \in \cH\ , \ \ k \in \bbN,  \label{G1}
\eea
form a singular sequence for the operator $\cS$, \eqref{4}, at the point 
$M^\flat = (1+ \lambda^\flat)^{-1}$. The first property of the Weyl criterion,
see, e.g., \cite[Thm.\,1,2]{BiSo}, 

\smallskip

\noindent
$1^\circ$. $\Vert U_k^\flat ; \cH \Vert =1$

\smallskip

\noindent
is just the normalization \eqref{G1}. The second condition 

\smallskip

\noindent
$2^\circ$. $U_k^\flat \to 0$ weakly in $\cH$
\noindent

\smallskip

\noindent
is not difficult either. Indeed, since the space $C_c^\infty( \overline{\Omega}
\smallsetminus \cO)$  of compactly supported infinitely smooth functions 
is dense in $\cH$, and, by definition,
\beas
\Vert \nabla u_k^\flat ; L^2(\Omega^\delta) \Vert
+ \Vert  u_k^\flat ; L^2(\partial \Omega^\delta \smallsetminus \Gamma^\delta ) 
\Vert \leq C_\delta
\eeas
for any fixed $\delta >0$, we conclude that for all $v \in C_c^\infty( 
\overline{\Omega} \smallsetminus \overline{\Pi}^\delta)$
\beas
( \nabla U_k^\flat , \nabla v )_\Omega + (U_k^\flat ,v)_{\partial \Omega} 
\to 0 
\eeas
because of the relation  $\Vert u_k^\flat ; L^2(\Omega) \Vert
= O(|\ln \varepsilon_k^\flat|^{1/2} )$, see \eqref{big}. So there remains to verify the property

\smallskip

\noindent
$3^\circ$. $\Vert \cS U_k^\flat  - M^\flat U_k^\flat ; \cH \Vert \to 0.$
\noindent

\smallskip

\noindent We repeat the calculation \eqref{Y5}, but since the supremum must now
be taken over the unit ball of $\cH$ instead of $\cH^\varepsilon$, we get
\beas
& & \Vert \cS U_k^\flat  - M^\flat U_k^\flat ; \cH \Vert 
= (1 + \lambda^\flat)^{-1} \Vert u_k^\flat; L^2(\Omega) \Vert^{-1} \nonumber \\
& & \times \sup \big| (\Delta u_k^\flat , V )_\Omega
- ( \partial_\nu u_k^\flat - \lambda^\flat u_k^\flat , V)_{\partial \Omega 
\smallsetminus \Gamma^\varepsilon} 
- (\partial_z u_k^\flat , V)_{\omega^\varepsilon}\big|.
\eeas
The first two terms inside the modulus have been estimated in Section \ref{sec4.5} 
by a bound approaching zero. To show that 
\bea
|(\partial_z u_k^\flat , V)_{\omega^\varepsilon} | \leq c \label{G2}
\eea
and thus to conclude with the proof of $3^\circ$, we need the following lemma
in addition to \ef{big}.

\BEL
\label{NiPo}
The trace inequality
\beas
\Vert V ; L^2(\omega^\varepsilon) \Vert &\leq &
c \sqrt{\varepsilon}\Vert V ; \cH  \Vert  
\eeas
holds  true for all $V \in \cH$, with constants $c$ depending on neither  $V \in \cH$ nor 
$\varepsilon$. 
\ENL

Proof. 
It is enough to prove the statement for smooth real-valued functions and by 
replacing $\Omega^\varepsilon \mapsto \Pi^d \smallsetminus \Pi^\varepsilon$ and
$\partial \Omega^\varepsilon \mapsto \Gamma^d \smallsetminus \Gamma^\varepsilon$. 
We use the coordinates $(\eta, z) = 
(z^{-2} y ,z )$ and the fundamental theorem of calculus 
\bea
& & \frac1\varepsilon \int\limits_{\omega^\varepsilon} |V(y,\varepsilon)|^2 dy
= \varepsilon^{-2(n-1)-1} 
\int\limits_{\omega} V(\varepsilon^2 \eta,\varepsilon)^2 d\eta
\roweq
\int\limits_{\omega}\big( z^{-2(n-1) -1}  V(z^2 \eta, z)^2 \big)
\Big|_{z = \varepsilon}  d \eta
\roweq
\int\limits_{\omega}\int\limits_\varepsilon^d
\frac{d}{dz} \Big( \chi_d(z) \big( z^{-2(n-1) -1}  V(z^2 \eta, z)^2 \big)
\Big) dz   d \eta
\rowleq
c \int\limits_\varepsilon^d \int\limits_{\omega}
z^{-2(n-1)}  \Big( (1 + z^ {-2} ) | V(z^2 \eta, z) | 
+ |\eta \cdot \nabla_y V(z^2 \eta ,z) |
\rowpl 
z^{-1} | \partial_z V(z^2 \eta,z) |  \Big) | V(z^2 \eta, z) | 
d\eta dz 
\rowleq
c \int\limits_\varepsilon^d \int\limits_{\omega^z}
\Big( (1 + z^ {-2} ) | V(y, z) |^2 
+ z^{-1} | V(y,z) | | \nabla V(y, z) | \Big)
dy dz  
\rowleq
c \int\limits_{\Pi^d \smallsetminus \Pi^\varepsilon}
\big( r^{-2}  | V(x) |^2 +  | \nabla V(x) |^2 \big) dx .
\label{O74}
\eea
The cut-off function $\chi_d \in C^\infty (\bbR)$ is fixed such that 
\beas
\chi_d(z) = 1 \ \mbox{for $z < d/3$ and $\chi_d(z) = 0$ for $ z > 2/3$.}
\eeas
Recalling \ef{32} at $\beta=0$, \ef{O74} yields the statement. \ \ 
$\boxtimes$

\bigskip

Now we obtain \ef{G2} from the following simple estimates, which are based on the 
properties of  the expressions contained in \eqref{Y5}:
\beas
& & \Vert \partial_z w_\pm ; L^2( \omega^\varepsilon)\Vert 
= c\varepsilon^{-n +1/2} ({\rm mes}_{n-1} \omega^\varepsilon)^{1/2}
= C \varepsilon^{-n+1/2} \varepsilon^{n-1} = C \varepsilon^{-1/2} , \\
& & \Vert \partial_z W_\pm ; L^2( \omega^\varepsilon)\Vert 
\leq  c \varepsilon^{-n+5/2} \varepsilon^{n-1} = c \varepsilon^{3/2}, \\
& & 
\Vert \varepsilon^{2-(n-3/2) \pm i \tau_0(\lambda^\flat)}
\partial_z Y_\pm ; L^2( \omega^\varepsilon)\Vert 
\leq  c \varepsilon^{-n+7/2} \varepsilon^{-2} \varepsilon^{n-1} = c 
\varepsilon^{1/2}, \\
& & \Vert \widetilde Z; L^2(\omega^\varepsilon)\Vert
\leq c \Vert r \partial_z^2 \widetilde Z; L^2(\Omega)\Vert \, 
\Vert \partial \widetilde Z; L^2(\Omega)\Vert  \leq C.
\eeas
The last estimate uses the trace inequality and the inclusions 
$\widetilde Z \in V_{-1}^1(\Omega)$ and $\nabla^2 \widetilde Z \in 
L_{1}^2(\Omega)$, proved in \cite{na549} and also mentioned in 
Theorem \ref{T2} and Remark \ref{R3}. 

Consequently, the Weyl criterion implies that any point $M^\flat
\in (0, M_\dagger)$ belongs to the essential spectrum of the operator
$\cS$, hence, $(\lambda_\dagger, + \infty) \subset \sigma_{\rm ess}$.
As has been outlined in Section \ref{sec1}, general results in \cite{Ko,MaPl2} 
imply  that the essential and continuous spectra of the Steklov problem coincide. 
Finally, the fact that the interval $(0, \lambda_\dagger)$ may only contain
points of the discrete  spectrum was shown in \cite{na405}. Thus, the above 
results  on blinking eigenvalues yield the formula $\sigma_{\rm co} = [\lambda_\dagger, + \infty)$,
which has already been obtained in \cite{na401}.

\begin{figure}
\label{fig2}
\begin{center}
{\includegraphics[height=6cm,width=10cm]{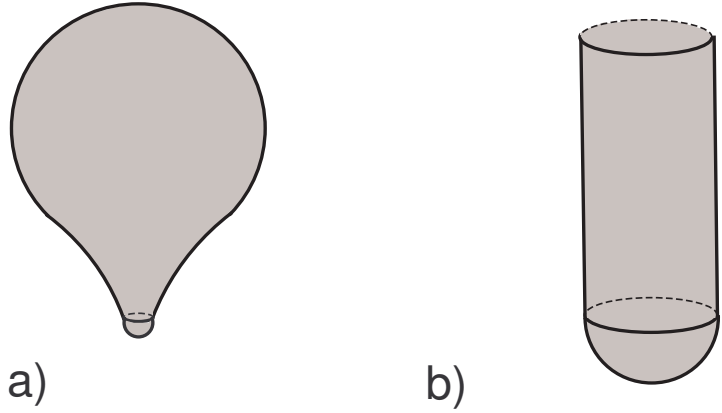}} 
\caption{a) Curved truncation surface, b) problem domain for the boundary 
layer. }
\end{center}
\end{figure}

\subsection{Other shapes of blunting.}
\label{sec5.2}
We consider the case where the truncation surface of the blunted cuspidal domain 
$\Omega^\varepsilon$ is defined by
\beas
\Upsilon^\varepsilon = \big\{ (y,z) \, : \, \big(\varepsilon^{-2} y ,
\varepsilon^{-2} (z- \varepsilon) \big) \in \Upsilon \big\}
\eeas
where $\Upsilon \subset \bbR_-^n$ is a piecewise smooth surface which touches
the half-space $\bbR_+^n$ at $\partial \omega \times \{ 0 \}$, see 
Fig.\,\ref{fig2}, a). Then, the spectrum of the Steklov-Dirichlet problem composed of 
the equations \eqref{8}, \ef{9} and 
\beas
u^\varepsilon(x) = 0, \, x \in \Upsilon^\varepsilon
\eeas
gets precisely the same properties as we established above for the domain
\eqref{7A} with the straight truncation surface. The only noteworthy  
modification  
in the proofs is related to the orthogonality condition \ef{28A} for the 
correction term $W$ in the asymptotic expansions near the cusp tip $\cO$;
these also appear in Section \ref{sec4.4}, where they provide the exponential
decay of the boundary layer terms $Y_\pm$. In the case of a curved truncation
surface as in Fig.\,\ref{fig2}, b), the boundary layer is to be found from the mixed boundary  value problem
\bea
& & - \Delta_\xi Y_\pm (\xi) = 0, \,  \xi \in \Xi_U, \ \ \ \ 
\partial_{\nu'} Y_\pm (\xi) = 0, \,  \xi \in \Sigma , 
\row
Y_\pm(\xi) = W_\pm(\xi',1), \,  \xi \in \Upsilon, 
\label{YY}
\eea
in the domain $\Xi_U$ which is bounded by the surfaces $\Sigma$ and $\Upsilon$
and contains the half-cylinder \eqref{Y3}. The homogeneous problem \ef{YY}
has a solution of the form
\beas
\cY(\xi) = |\omega|^{-1} \xi_n + C_\cY + O( e^{-\beta_\Xi \xi_n})
\ \ \mbox{as} \ \xi_n \to + \infty. 
\eeas
Then, the exponential decay of the solution $Y_\pm$ of \ef{YY}
is supported by the orthogonality condition 
\bea
\int\limits_\Upsilon W_\pm(\xi',1) \partial_\nu \cY(\xi) ds_\xi = 0  ,
\label{newY}
\eea
which replaces \ef{28A} everywhere. Note that in the case of a straight end
$\omega^\varepsilon$ we have $\cY (\xi) = |\omega|^{-1} \xi_n$ so that 
formula \eqref{newY} turns into \eqref{29A}.

\subsection{Other boundary conditions at the end $\omega^\varepsilon$. }
\label{sec5.3}
If the Dirichlet condition \ef{10} is replaced by the Neumann condition
\bea
\partial_z u^\varepsilon (x) = 0 \ , \ \ x\in \omega^\varepsilon, \label{Ne}
\eea
the phenomena of blinking and gliding are preserved, but the relationship
\ef{62} is changed a bit, since the normal derivative of \ef{61B} vanishes
at $z = \varepsilon$ provided
\beas
- 2 \tau_0(\lambda^\flat) \ln \varepsilon = \Theta(\lambda^\flat)
+ \vartheta (\lambda^\flat) + \pi 
\ \ \ ( {\rm mod} \, 2 \pi) ,
\eeas
where 
\beas
e^{i \vartheta_n (\lambda^\flat)} = 
\frac{(n - 3/2) - i \tau_0(\lambda^\flat)}{(n - 3/2) + i \tau_0(\lambda^\flat)}.
\eeas
The same formula and similar conclusions occur in the case of the Steklov
condition 
\beas
\partial_z u^\varepsilon(x) = \lambda^\varepsilon u^\varepsilon(x) , 
\ \ x  \in \omega^\varepsilon ;
\eeas
the proofs require some minor modifications.

\end{document}